\documentclass[10pt]{amsart} 
\title{Tableaux combinatorics for the asymmetric exclusion process}
\author{Sylvie Corteel and Lauren K. Williams}

\address{CNRS LRI, Universit\'e Paris-Sud, B\^atiment 490, 91405 Orsay
Cedex France} 
\email{Sylvie.Corteel@lri.fr}
\address{University of California, Berkeley, and Harvard University, Cambridge}
\email{lauren@math.harvard.edu}



\usepackage{pstricks, pst-node, amssymb}

\psset{unit=1pt, arrowsize=4pt, linewidth=.7pt}
\psset{linecolor=blue}
\newgray{grayish}{.90}
\newrgbcolor{embgreen}{0 .5 0}
\def\Le{\hbox{\rotatedown{$\Gamma$}}}
\def\vblack(#1, #2)#3{\cnode*[linecolor=black](#1, #2){3}{#3}}
\def\vwhite(#1,#2)#3{\cnode[linecolor=black,fillcolor=white,fillstyle=solid](#1,
#2){3}{#3}}
\countdef\x=23
\countdef\y=24
\countdef\z=25
\countdef\t=26

\def\tbox(#1,#2)#3{
\x=#1 \y=#2
\multiply\x by 12
\multiply\y by 12
\z=\x \t=\y
\advance\z by 12
\advance\t by 12
\psline(\x,\y)(\x,\t)(\z,\t)(\z,\y)(\x,\y)
\advance\x by 6
\advance\y by 6
\rput(\x,\y){{\bf #3}}}

\usepackage{amsmath, amsthm, amssymb, amsbsy}
\usepackage{amsfonts, latexsym, stmaryrd, amscd, xy}
\usepackage[mathscr]{eucal}
\usepackage{epsfig}

\newtheorem{theorem}{Theorem}[section]

\newtheorem{proposition}[theorem]{Proposition}
\newtheorem{lemma}[theorem]{Lemma}

\newtheorem{example}[theorem]{Example}
\newtheorem{corollary}[theorem]{Corollary}

\newtheorem{remark}[theorem]{Remark}

\newtheorem{definition}[theorem]{Definition}

\usepackage{amsfonts}
\usepackage{xy}
\usepackage{amssymb}
\usepackage{amsmath}
\newcommand{\T}{\mathcal{T}}
\newcommand{\ttt}{\tau}

\DeclareMathOperator{\expanse}{expanse}
\DeclareMathOperator{\wt}{wt}

\newcommand{\thmrefer}[1]{\renewcommand\thetheorem
  {\protect\ref{#1}}\addtocounter{theorem}{-1}}

\xyoption{all}
\CompileMatrices

\begin{document}



\begin{abstract}

The partially asymmetric exclusion process (PASEP) is an important model
from statistical mechanics which describes a system of interacting 
particles hopping left and right on a one-dimensional lattice of $n$
sites.  It is partially asymmetric in the sense that the probability 
of hopping left is $q$ times the probability of hopping right.  
Additionally, particles may enter from the left with probability 
$\alpha$ and exit from the right with probability $\beta$.

In this paper we prove a close connection between the PASEP and the 
combinatorics of permutation tableaux.  (These tableaux come indirectly 
from the totally nonnegative part of the Grassmannian, via work of 
Postnikov, and were studied in a paper of Steingrimsson 
and the second author.)  Namely, we prove that in the long time 
limit, the probability that the PASEP is in a particular configuration 
$\tau$ is essentially the generating function for permutation tableaux
of shape $\lambda(\tau)$ enumerated according to three statistics.
The proof of this result uses a result of 
Derrida, Evans, Hakim, and Pasquier on the {\it matrix ansatz} for the 
PASEP model.  

As an application, we prove some monotonicity results for the 
PASEP.  We also derive some enumerative
consequences for permutations enumerated according to various 
statistics such as weak excedence set, descent set, crossings,
and occurences of generalized patterns.
\end{abstract}

\maketitle

\section{Introduction}

The partially asymmetric exclusion process (PASEP) 
is an important model from statistical mechanics which is quite  simple 
but surprisingly rich: it exhibits 
boundary-induced phase transitions, spontaneous symmetry 
breaking, and phase separation.  The PASEP is regarded as a primitive
model for biopolymerization \cite{bio}, traffic flow \cite{traffic}, 
and formation of shocks \cite{shock};
it also appears in a kind of sequence alignment problem in 
computation biology \cite{comput}.
More recently it has been noticed that the PASEP model has relations
to orthogonal polynomials \cite{Sasamoto}, and to interesting 
combinatorial phenomena \cite{jumping, Corteel, Zeilberger}.
The goal of this 
paper is to prove a precise 
connection betweeen the PASEP model and 
permutation tableaux, certain $0-1$ tableaux introduced in \cite{SW}.

The PASEP model describes a system of particles hopping left and right on a one-dimensional lattice of $n$ sites.  Particles may enter the system from 
the left with a rate $\alpha dt$ and may exit the system from the right 
at a rate $\beta dt$.  The probability of hopping left is $q$ times the 
probability of hopping right.

Let $\frac{f_{\tau}(q)}{Z_n}$ denote the probability that the PASEP model is in 
a particular configuration $\tau$ in the steady state.  Here, $Z_n$ is the 
{\it partition function} of the PASEP model.
In this paper
we consider a solution $(D_1,E_1,V_1,W_1)$ 
to the ``matrix ansatz" for the PASEP model which naturally 
relates to permutation tableaux.  That is, 
we find matrices $D_1$ and $E_1$ and vectors
$V_1$ and $W_1$
which satisfy the relations of Theorem \ref{ansatz}; we then
prove that expressions of the form 
$W_1 (\prod_{i=1}^n (\ttt_i D_1 + (1-\ttt_i)E_1))V_1$,
where $\ttt=(\ttt_1, \dots , \ttt_n)\in \{0,1\}^n$, 
are (Laurent) polynomials in $q, \alpha, \beta$ which
enumerate permutation tableaux of a fixed shape 
$\lambda(\tau)$ according to three statistics.
Using a result of Derrida et al \cite{Derrida1}, we are
then able to prove our main result: that  
$f_{\tau}(q)$ 
is the 
generating function for permutation tableaux of a fixed shape
$\lambda(\tau)$.  

It then follows from work of 
the second author and Steingr\'{\i}msson \cite{SW} that $f_{\tau}(q)$ is 
also the generating function for:
permutations in $S_{n+1}$ with 
a fixed set $W(\tau)$ of weak excedences, enumerated according to crossings;
permutations in $S_{n+1}$ with a fixed set 
$D(\tau)$ of descents, enumerated according to occurrences of the 
{\it generalized pattern} $2-31$.
Additionally, these results imply that 
the expression 
$W_1 (D_1+E_1)^n V_1$ for the partition function $Z_n$ is also 
the weight-generating function for 
permutations in $S_{n+1}$,
enumerated according to crossings, and 
permutations in $S_{n+1}$,
enumerated according to occurrences of the 
{\it generalized pattern} $2-31$.

Our main result refines the theorem of the first author \cite{Corteel}, who
showed that if $\alpha = \beta = 1$, then
in the steady state, the probability that the model
is in a configuration with $k$ occupied sites is equal to 
$\frac{\hat{E}_{k+1,n+1}(q)}{Z_n}$, where 
$\hat{E}_{k+1,n+1}(q)$ is the $q$-Eulerian polynomial introduced
by the second author \cite{Williams}.

The structure of this paper is as follows.  
In Section \ref{setup} we define the PASEP model and review
the ``matrix ansatz," presenting a solution
$(D_1, E_1)$ which has an interpretation in terms of {\it permutation 
tableaux.}
In Section \ref{PermTableaux}
we define  permutation tableaux, certain $0-1$ tableaux which 
were defined in \cite{SW} and which are naturally in bijection
with permutations.  We then 
show how the solution $(D_1,E_1)$ to the matrix ansatz leads to 
a natural connection between the PASEP model and permutation tableaux,
and hence to permutations.  In Section \ref{Proof}, we prove 
the main result of the previous section, and in Section 
\ref{Applications}, we give some applications.
Finally, in section \ref{Motzkin}, we show how a classical solution 
$(D_0,E_0)$ of the matrix ansatz leads to a natural connection 
between the PASEP model and bicolored Motzkin paths.
This recovers results of Brak et al \cite{Corteel2}.

It is interesting to note that permutation tableaux are closely 
connected to total positivity for the Grassmannian
\cite{Williams, Postnikov}.  This suggests an intriguing 
connection between total positivity and the PASEP model.

\textsc{Acknowledgments:}
We would like to thank the referee for insightful comments.
The second author is grateful to Persi Diaconis for interesting discussions
in Banff 
and for encouraging her to work on this project.
She would also like to thank Kolya Reshetikhin
and Richard Stanley for encouragement and useful comments.

\section{The PASEP model and the ``matrix ansatz"}\label{setup}

In the physics literature, the PASEP is defined as follows.
\begin{definition}
We are given a one-dimensional
lattice of $N$ sites, such that each site $i$ ($1 \leq i \leq N)$
is either occupied by a particle ($\tau_i=1$) or is empty
($\tau_i=0$).  At most one particle may occupy a given site.
During each infinitesimal time interval $dt$, each particle
in the system has a probability $dt$ of jumping to the next
site on its right (for particles on sites $1 \leq i \leq N-1$) and
a probability $q dt$ of jumping to the next site on its left
(for particles on sites $2 \leq i \leq N$).  Furthermore, a particle
is added at site $i=1$ with probability $\alpha dt$ if site $1$  
is empty and a particle is removed from site $N$ with probability
$\beta dt$ if this site is occupied.
\end{definition}

Note that we will sometimes denote a state of the PASEP as a
word in $\{0,1\}^N$ and sometimes as a word in $\{\circ,\bullet\}^N$.
In the latter notation, the symbol $\circ$ denotes the absence of
a particle, which one can also think of as a white particle.

It is not too hard to see \cite{jumping} that our previous formulation
of the PASEP is equivalent to the following
discrete-time Markov chain.

\begin{definition}
Let $\alpha$, $\beta$, and $q$ be constants such that 
$0 < \alpha \leq 1$, $0 < \beta \leq 1$, and $0 \leq q \leq 1$.
Let $B_N$ be the set of all $2^N$ words in the
language $\{\circ, \bullet\}^*$.
The PASEP is the Markov chain on $B_N$ with
transition probabilities:
\begin{itemize}
\item  If $X = A\bullet \circ B$ and
$Y = A \circ \bullet B$ then
$P_{X,Y} = \frac{1}{N+1}$ (particle hops right) and
$P_{Y,X} = \frac{q}{N+1}$ (particle hops left).
\item  If $X = \circ B$ and $Y = \bullet B$
then $P_{X,Y} = \frac{\alpha}{N+1}$ (particle enters from left).
\item  If $X = B \bullet$ and $Y = B \circ$
then $P_{X,Y} = \frac{\beta}{N+1}$ (particle exits to the right).
\item  Otherwise $P_{X,Y} = 0$ for $Y \neq X$
and $P_{X,X} = 1 - \sum_{X \neq Y} P_{X,Y}$.
\end{itemize}
\end{definition}

\begin{figure}[h]
\centering
\includegraphics[height=1.3in]{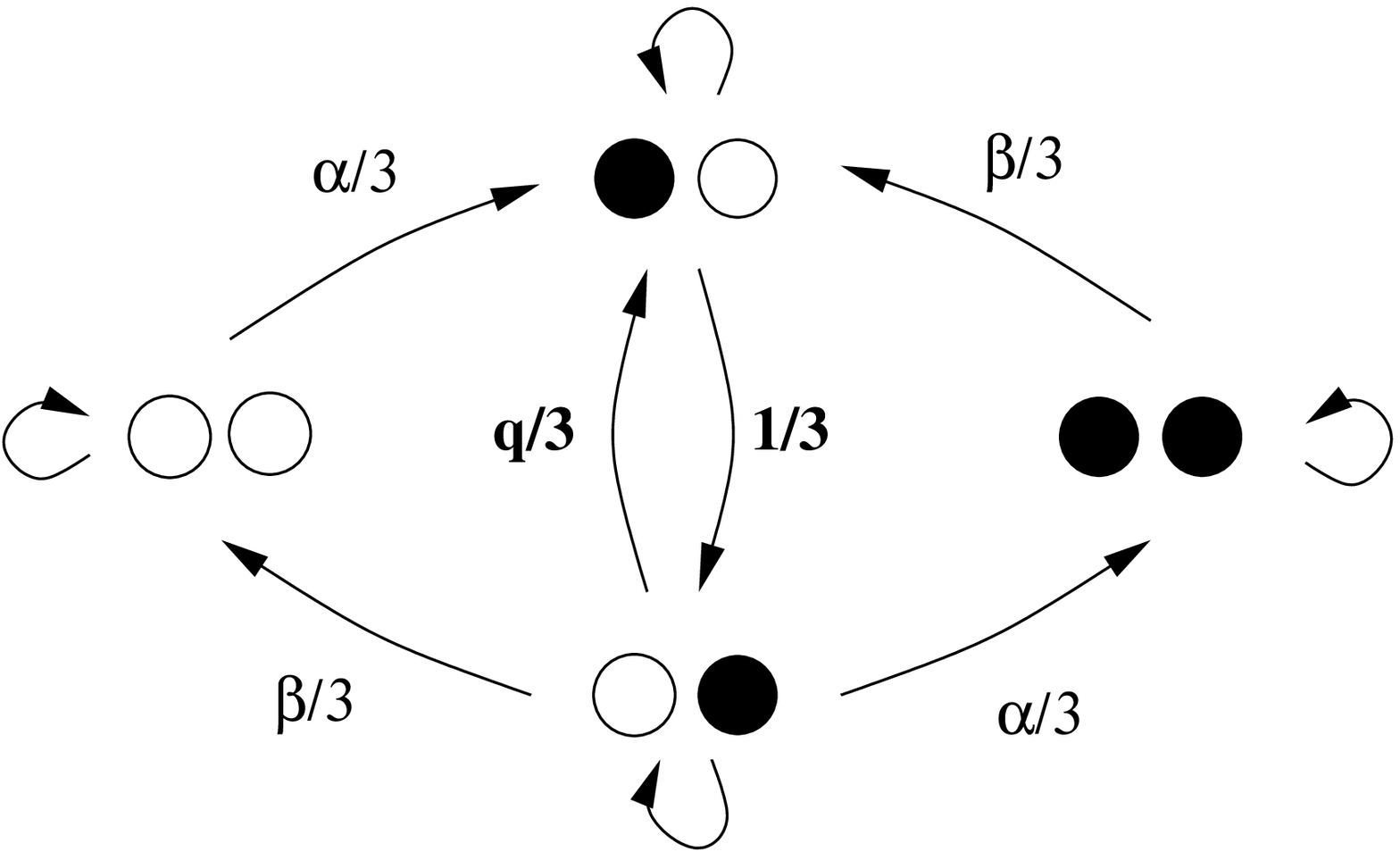}
\caption{The state diagram  of the PASEP model for $n=2$}
\label{states}
\end{figure}

See Figure \ref{states} for an
illustration of the four states, with transition probabilities,
for the case $n=2$.  Note that the probabilities of the loops
are determined from the figure by the fact that the sum of the probabilities
on all outgoing arrows from a given state must be $1$.  So for example
the probability on the bottom-most loop is 
$1-\frac{q}{3}-\frac{\alpha}{3}-\frac{\beta}{3} = \frac{3-q-\alpha-\beta}{3}$.

In the long time limit, the system reaches a steady state where all 
the probabilities $P_n(\ttt_1, \ttt_2, \dots , \ttt_n)$ of finding
the system in configurations $(\ttt_1, \ttt_2, \dots , \ttt_n)$ are
stationary, i.e.\ satisfy 
\begin{equation*}
\frac{d}{dt} P_n(\tau_1, \dots , \ttt_n) = 0.
\end{equation*}
Moreover, the stationary distribution is unique \cite{Derrida1}, as shown
by Derrida et al.

The question is now to solve for the probabilities 
$P_n(\ttt_1, \dots , \ttt_n)$.  For convenience, we 
define unnormalized weights $f_n(\ttt_1, \dots , \ttt_n)$, which are
equal to the $P_n(\ttt_1, \dots , \ttt_n)$ up to a constant:
\begin{equation*}
P_n(\ttt_1, \dots , \ttt_n) = f_n(\ttt_1, \dots , \ttt_n)/Z_n,
\end{equation*}
where 
$Z_n$ is the {\it partition function}
$\sum_{\tau} f_n(\ttt_1, \dots , \ttt_n)$.  The sum defining
$Z_n$ is 
over all possible configurations $\ttt \in \{0,1\}^n$.

The ``matrix ansatz" has been used by Derrida et al \cite{Derrida1}
to obtain exact expressions for all the $P_n(\ttt_1, \dots , \ttt_n)$.

More precisely, they show the following.

\begin{theorem} \cite{Derrida1} \label{ansatz}
Suppose that $D$ and $E$ are matrices,  $V$ is a column vector,
and $W$ is a row vector,
such that the following conditions hold:
\begin{align*}
& DE - qED = D+E \\
& DV = \frac{1}{\beta} V \\
& WE = \frac{1}{\alpha} W
\end{align*}
Then 
\begin{equation*}
f_n(\ttt_1, \dots , \ttt_n) = W (\prod_{i=1}^n (\ttt_i D + (1-\ttt_i)E))V.
\end{equation*}
\end{theorem}

Note that $\prod_{i=1}^n (\ttt_i D + (1-\ttt_i)E)$ is simply
a product of $n$ matrices $D$ or $E$ with matrix $D$ at position $i$
if site $i$ is empty ($\ttt_i=0)$.

\begin{remark}\label{partition}
It follows from Theorem \ref{ansatz} that the partition function
$Z_n$ is equal to $W (D+E)^n V$.
\end{remark}

We will now describe a solution $(D_1,E_1,V_1,W_1)$
to the matrix ansatz. 
It seems to have not been considered before, although it naturally
generalizes a solution given in \cite{Derrida1} for the 
$q=0$ case.  Our solution
has an interpretation in 
terms of permutation tableaux and we found it through consideration
of the combinatorics of these tableaux; however, one could 
also arrive at these matrices by choosing the basis 
$\{V_1, EV_1, E^2 V_1, \dots \}$ for the infinite-dimensional vector
space that $D$ and $E$ act on, and then writing $D$ and $E$ 
according to this basis.  (This method actually
yields the transposes of the matrices
$D_1$ and $E_1$ which we use below.)

Let $D_1$ be the (infinite) upper triangular
matrix $(d_{ij})$ such that 
$d_{i,i+1}=\beta^{-1}$ and $d_{i,j}=0$ for $j\neq i+1$.

That is, $D_1$ is the matrix
\[ \left( \begin{array}{ccccc}
0 & \beta^{-1} & 0 & 0 & \dots \\
0 & 0 & \beta^{-1} & 0 & \dots \\
0 & 0 & 0 & \beta^{-1} & \dots \\
0 & 0 & 0 & 0 & \dots \\
\vdots & \vdots & \vdots & \vdots &
\end{array} \right). \]

Let $E_1$ be the (infinite) lower triangular matrix
$(e_{ij})$ such that for $j \leq i$, 
$e_{ij} = \beta^{i-j} (\alpha^{-1} q^{j-1} {i-1 \choose j-1} +
\sum_{r=0}^{j-2} {i-j+r \choose r} q^r)$.
Otherwise, $e_{ij} = 0$.

That is, $E_1$ is the matrix 
\[ \left( \begin{array}{ccccc}
\alpha^{-1} & 0 & 0 & 0 &  \dots \\
\alpha^{-1} \beta & 1+\alpha^{-1} q & 0 & 0 & \dots \\
\alpha^{-1} \beta^2 & \beta(1+2\alpha^{-1} q) & 1+q+\alpha^{-1} q^2 & 0&
\dots \\
\alpha^{-1} \beta^3  & \beta^2 (1+3\alpha^{-1}q) &
\beta(1+2q+3\alpha^{-1} q^2)  & 1+q+q^2+\alpha^{-1} q^3 & \dots \\
\vdots & \vdots & \vdots & \vdots &
\end{array} \right). \]

Observe that when $\alpha = \beta = 1$, we have 
$e_{ij} = 
 \frac{[i]^{(i-j)}}{(i-j)!}$.
Here, $[i]^{(k)}$ represents the $k$th derivative
of $[i]$ with respect to $q$, and $[i]$ is the $q$-analog of 
the number $i$, namely $1+q+\dots + q^{i-1}$.

And then $E_1$ becomes the matrix 
\[ \left( \begin{array}{cccccc}
1 & 0 & 0 & 0 & 0& \dots \\
1 & [2] & 0 & 0 & 0&\dots \\
1 & [3]' & [3] & 0 & 0&\dots \\
1 & \frac{[4]''}{2} & [4]' & [4] & 0&\dots \\
1 & \frac{[5]'''}{6} & \frac{[5]''}{2} & [5]' & [5] & \dots \\
\vdots & \vdots & \vdots & \vdots & \vdots &
\end{array} \right). \]

Let $W_1$ be the (row) vector $(1,0,0,\dots )$ and $V_1$ be the (column)
vector
$(1,1,1,\dots)$.
It is now easy to check that the required relations hold.

\begin{lemma}
With the definitions of $D_1, E_1, V_1, W_1$ above, the following
relations hold:
$D_1 E_1 -q E_1 D_1 = D_1 + E_1$, 
$DV_1 = \frac{1}{\beta} V_1$, and $W_1 E = \frac{1}{\alpha} W_1$.
\end{lemma}

\begin{proof}
First note that 
\begin{equation*} (D_1 E_1)_{i,j}=\beta^{-1} (E_1)_{i+1,j} =
 \beta^{i-j} \left(\alpha^{-1} q^{j-1} {i \choose j-1} +
\sum_{r=0}^{j-2} {i-j+r+1 \choose r} q^r \right),
\end{equation*}
when $j \leq i+1$, and is equal to $0$ otherwise.
Then note that 
\begin{equation*}
q(E_1 D_1)_{i,j} = q \beta^{-1} (E_1)_{i,j-1} = 
 \beta^{i-j} q \left(\alpha^{-1} q^{j-2} {i-1 \choose j-2} +
\sum_{r=0}^{j-3} {i-j+r+1 \choose r} q^r\right).
\end{equation*}
when $1 \leq j-1 \leq i$, and is equal to $0$ otherwise.

Putting these together, we find that 
$(D_1 E_1 - qE_1 D_1)_{i,j}$ is equal to  
\begin{align*}
 &= 
 \beta^{i-j} \left(\alpha^{-1} q^{j-1} {i \choose j-1} +
\sum_{r=0}^{j-2} {i-j+r+1 \choose r} q^r\right) - \\
  & \qquad \qquad \qquad \qquad \qquad  \beta^{i-j} \left(\alpha^{-1} q^{j-1} {i-1 \choose j-2} +
   \sum_{r=0}^{j-3} {i-j+r+1 \choose r} q^{r+1}\right) \\
   &= 
\beta^{i-j} \alpha^{-1} q^{j-1} {i-1 \choose j-1} 
+ \beta^{i-j} \sum_{r=0}^{j-2} {i-j+r+1 \choose r} q^r - 
\beta^{i-j} \sum_{s=1}^{j-2} {i-j+s \choose s-1} q^s \\ 
   &= 
\beta^{i-j} \alpha^{-1} q^{j-1} {i-1 \choose j-1} 
+ \beta^{i-j} \sum_{r=0}^{j-2} {i-j+r \choose r} q^r,
\end{align*}
when $1 < j \leq i$.
Also, we have that 
$(D_1 E_1 - qE_1 D_1)_{i,j}=\beta^{-1}$ 
when $j=i+1$, is equal to $\alpha^{-1}\beta^{i-1}$ when  $j=1$, 
and is equal to 
$0$ when $j>i+1$.
These are precisely the matrix entries of 
$D_1+E_1$.
\end{proof}

As we will show in Sections \ref{PermTableaux} and \ref{Proof}  
the matrix product $W_1 (\prod_{i=1}^n (\ttt_i D_1 + (1-\ttt_i)E_1))V_1$
from Theorem \ref{ansatz}
has a combinatorial interpretation as a generating function
for  permutation tableaux.

\begin{example}\label{example1}
Suppose that we are considering the PASEP model with $3$ sites, and are
interested in computing the probability that in the long time limit,
the system reaches the configuration which has a black particle in 
the second  site (but the first and third sites are empty). 
This configuration is represented by the vector 
$\tau = (0,1,0)$.  
Using Theorem \ref{ansatz}, we see that $f_3(0,1,0)$ is
equal to the matrix
product $W_1 (E_1 D_1 E_1) V_1$;
this product is $\alpha^{-2} + \alpha^{-1}\beta^{-1} + \alpha^{-2}\beta^{-1}q$.
And by Remark \ref{partition}, the 
partition function $Z_3$ is equal to 
$W_1 (D_1 +E_1)^3 V_1$, which is in this case
$\alpha^{-3}+2\alpha^{-2}+2\alpha^{-1}+\alpha^{-2}\beta^{-1}+
2\alpha^{-1}\beta^{-1} +2\beta^{-1}+\alpha^{-1}\beta^{-2} + 2\beta^{-2}+
\beta^{-3} +
q(\alpha^{-2}+\alpha^{-2}\beta^{-1}+4\alpha^{-1}\beta^{-1}+\alpha^{-1}\beta^{-2}+\beta^{-2}) +
q^2(\alpha^{-2}\beta^{-1} + \alpha^{-1} \beta^{-2}).$
Therefore
the probability that in the long time limit the system is in 
the configuration $(0,1,0)$ is 
$\frac{\alpha^{-2} + \alpha^{-1}\beta^{-1} + \alpha^{-2}\beta^{-1}q}{Z_3}$.
\end{example}

\section{Connection with permutation tableaux}\label{PermTableaux}

Recall that a {\em partition} $\lambda = (\lambda_1, \dots,
\lambda_k)$ is a weakly decreasing sequence of nonnegative
integers. For a partition $\lambda$, where $\sum \lambda_i = m$, the
{\em Young diagram} $Y_\lambda$ of shape $\lambda$ is a left-justified
diagram of $m$ boxes, with $\lambda_i$ boxes in the $i$-th row.
We define the {\it expanse} of $\lambda$ or $Y_\lambda$ to be the sum 
of the number of rows and the number of columns.  Note that 
we will allow a row to have length $0$-- i.e.\ we allow $\lambda_i=0$ --
and we distinguish two partitions that differ in their number of
empty rows.

\begin{figure}[h]
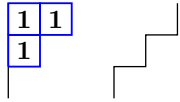

\pspicture(0,0)(40,40)        
\psline[linecolor=black,linewidth=0.5pt]{-}(0,0)(0,36)
\tbox(0,2){1}                 
\tbox(0,1){1}                 
\tbox(1,2){1}                 
\endpspicture
\pspicture(0,0)(40,40)
\psline[linecolor=black,linewidth=0.5pt]{-}(0,0)(0,12)(12,12)(12,24)(24,24)(24,36)
\endpspicture
\caption{Young diagram and path for $\lambda=(2,1,0)$}
\label{path}
\end{figure}

We will often identify
a Young diagram $Y_\lambda$ with expanse  $n$ 
with the lattice path $p(\lambda)$ of length $n$
which takes unit steps south and west,
beginning at
the north-east corner of $Y_\lambda$ and ending at
the south-west corner.
Note that such a lattice path always begins with a step south.
See Figure \ref{path} for the path corresponding to the 
Young diagram of shape $(2,1,0)$.

If $\tau \in \{0,1\}^{n-1}$, then we associate to it a Young diagram 
$\lambda(\tau)$ with expanse $n$ as follows.  First we define 
a path $p = (p_1, \dots , p_n) \in \{S,W\}^n$ such that $p_1=S$, and 
$p_{i+1}=S$ if and only if $\tau_i = 1$.  We then define
$\lambda(\tau)$ to be the partition associated to this path $p$.
We denote by $\phi$ the inverse of the above map: 
it is a bijection  
from the set of Young diagrams
with expanse $n$ to the set of $n-1$-tuples in $\{0,1\}^{n-1}$.

As in \cite{SW}, we
define a {\em permutation tableau} $\T$ to be a partition $\lambda$
together with a filling of the boxes of $Y_\lambda$ with $0$'s and
$1$'s such that the following properties hold:
\begin{enumerate}
\item Each column of the rectangle contains at least one $1$.
\item There is no $0$ which has a $1$ above it in the same column
{\em and} a $1$ to its left in the same row.
\end{enumerate}

We call such a filling a {\em valid} filling of $Y_\lambda$.  

Note that the second requirement above can be rephrased in the following
way.  Read the columns of a permutation tableau $\T$ from right to left.  
If in any column we have a $0$ which lies beneath some $1$, then all 
entries to the left of $0$ (which are in the same row) must also
be $0$'s.

Note that if we forget the requirement  (1) above we recover the definition of
a $\Le$-diagram~\cite{Postnikov, Williams}, 
an object which represents
a cell in the totally nonnegative part of the Grassmannian.

We will now define a few statistics on permutation tableaux.
We define the {\it weight} $\wt(\T)$ of a permutation tableau $\T$ with k 
columns to be the total number of $1$'s in the filling minus $k$.
(We subtract $k$ since there must be at least $k$ $1$'s in 
a valid filling of a tableau with $k$ columns.)
In other words, we are counting how many extra $1$'s each column contains
beyond the requisite one.

We define $f(\T)$ to be the number of $1$'s in the first
row of $\T$.

We say that an entry in a column of a permutation tableau is 
{\it restricted} if 
that entry is a $0$ which lies below some $1$. 
And we say that a row is 
{\it unrestricted}
if it does not contain a restricted entry.  
Define $u(\T)$ to be the number of unrestricted rows of $\T$ 
minus $1$.  (We subtract $1$ since the top row of a tableau
is always unrestricted.)

Figure \ref{PermTab} gives an example of a permutation tableau $\T$
with weight $19-10=9$ and expanse $17$, such that 
$u(\T) = 3$ and $f(\T) = 5$.

\begin{figure}[h]
\pspicture(-95,-12)(230,94)
\rput(190,36)
{$\begin{array}{l}
\lambda=(10,9,9,8,5,2)
\end{array}$}
\psline[linecolor=black,linewidth=0.5pt]{-}(0,-12)(120,-12)(120,72)(0,72)(0,-12)
\tbox(0,0){1}
\tbox(1,0){1}
\tbox(0,1){0}
\tbox(1,1){0}
\tbox(2,1){0}
\tbox(3,1){1}
\tbox(4,1){1}
\tbox(0,2){0}
\tbox(1,2){0}
\tbox(2,2){0}
\tbox(3,2){0}
\tbox(4,2){0}
\tbox(5,2){0}
\tbox(6,2){1}
\tbox(7,2){1}
\tbox(0,3){0}
\tbox(1,3){0}
\tbox(2,3){0}
\tbox(3,3){0}
\tbox(4,3){0}
\tbox(5,3){0}
\tbox(6,3){0}
\tbox(7,3){0}
\tbox(8,3){0}
\tbox(0,4){1}
\tbox(1,4){1}
\tbox(2,4){1}
\tbox(3,4){1}
\tbox(4,4){0}
\tbox(5,4){1}
\tbox(6,4){1}
\tbox(7,4){1}
\tbox(8,4){1}
\tbox(0,5){0}
\tbox(1,5){1}
\tbox(2,5){1}
\tbox(3,5){0}
\tbox(4,5){0}
\tbox(5,5){1}
\tbox(6,5){0}
\tbox(7,5){1}
\tbox(8,5){0}
\tbox(9,5){1}
\endpspicture
\caption{A permutation tableau}
\label{PermTab}
\end{figure}

Now we will be interested in the question of enumerating permutation 
tableaux according to their shape, weight, unrestricted rows, and 
first row.
That is, we are interested in computing the polynomials
$F_{\lambda}(q):= \sum_{\T} q^{\wt(\T)} \alpha^{-f(\T)} \beta^{-u(\T)}$, 
where the 
sum ranges over all permutation tableaux $\T$ of shape
$\lambda$.  Let us also define the 
polynomials $F^n(q):= 
\sum_{\T} q^{\wt(\T)} \alpha^{-f(\T)} \beta^{-u(\T)}$, 
where the sum ranges over all permutation tableaux $\T$ with expanse  $n$.

Our main result is the following.

\begin{theorem}\label{Th1}
Fix a partition $\lambda$ with expanse $n+1$.
Let $(\ttt_1, \dots , \ttt_n) \in \{0,1\}^n$ be 
$\phi(\lambda)$.  Then 
\begin{equation*}
F_{\lambda}(q) = W_1 (\prod_{i=1}^n (\ttt_i D_1 + (1-\ttt_i) E_1)) V_1. 
\end{equation*}
Moreover, the generating function $F^{n+1}(q)$ for {\it all}
permutation tableaux with expanse $n+1$ is 
$W_1 (D_1+E_1)^n V_1$.
\end{theorem}

An examination of the formula in Theorem \ref{Th1} reveals that
$F_{\lambda}(q)$ is just the sum of the entries
in the top row of a certain product of $n$ matrices, each one either
$D_1$ or $E_1$, where this product corresponds via $\phi$ to the 
shape of the partition $\lambda$.

We remark that Theorem \ref{Th1} is particularly 
nice because it is a {\it positive}
formula.  In \cite{SW}, a recurrence was given for 
$F_{\lambda}(q)$ (in the case where $\alpha=\beta=1$), 
but not an explicit formula.

\begin{example}\label{example2}
We now illustrate Theorem \ref{Th1} with an example.
Suppose that we want to calculate the weight generating function
for permutation tableaux of shape $\lambda = (2,1)$.  Note that
$\phi(\lambda)=(0,1,0)$.  Therefore we need to evaluate
the expression $W_1 E_1 D_1 E_1 V_1$.  This is equal to 
$\alpha^{-2}\beta^{-1}q +
\alpha^{-2}+\alpha^{-1}\beta^{-1}$.  
And indeed, as shown in Figure \ref{Hello}, there are three
permutation tableaux of shape $(2,1)$, whose statistics correspond to 
the three terms above.
Compare our results here to 
Example
 \ref{example1}.

\begin{figure}[h]
\pspicture(0,0)(40,40)        
\tbox(0,2){1}                 
\tbox(0,1){1}                 
\tbox(1,2){1}                 
\endpspicture
\pspicture(0,0)(40,40)
\tbox(0,2){1}
\tbox(1,2){1}
\tbox(0,1){0}
\endpspicture
\pspicture(0,0)(40,40)
\tbox(0,2){0}
\tbox(1,2){1}
\tbox(0,1){1}
\endpspicture

\caption{Permutation tableaux of shape $(2,1,0)$}
\label{Hello}
\end{figure}
\end{example}

We will defer the proof of Theorem \ref{Th1} to Section \ref{Proof},
in order to first explain various consequences of the result.

Theorem \ref{Th1}  
together with Theorem \ref{ansatz} implies the following.

\begin{corollary}\label{Cor1}
Fix $\ttt = (\ttt_1, \dots , \ttt_n) \in \{0,1\}^n$, and 
let $\lambda := \lambda(\tau)$.  (Note that $\expanse (\lambda)=n+1$.)
The probability of finding the PASEP model in configuration 
$(\ttt_1, \dots , \ttt_n)$ in the steady state is 
\begin{equation*}
\frac{F_{\lambda}(q)}{Z_n}.
\end{equation*}
Here, $F_{\lambda}(q)$ is the weight-generating function for 
permutation tableaux of shape $\lambda$. 
Moreover, the partition function $Z_n$ for the PASEP model
is equal to the generating function for all permutation
tableaux with expanse $n+1$, enumerated according to 
weight, unrestricted rows, and the first row.
\end{corollary}

In the case that $\alpha=\beta=1$, 
this corollary together with results from work  
with
Steingr\'{\i}msson \cite{SW} 
will give us two more combinatorial interpretations
of the probabilities in the steady state, in terms of 
excedences and crossings, and descents and patterns in permutations.  
We now recall 
some definitions and results.

If $\pi \in S_n$ is a permutation, we say that $\pi$ has a 
{\it weak excedence} in position $i$ if $\pi(i) \geq i$.  
Note that there is always a weak excedence in position $1$.
The {\it weak excedence set} of $\pi$ is the set of all indices
$i$ where there is a weak excedence.
Additionally, we say that $\pi$ has a {\it crossing} in positions  $(i,j)$ 
if either $j < i \leq \pi(j) < \pi(i)$ or 
$\pi(i) < \pi(j) < i < j$.
Crossings were first defined by the first author in \cite{Corteel}.

\begin{theorem} \cite[Theorem 7]{SW} \label{Bij}
Let $T(k,n,c)$ be the set of permutation tableaux with 
$k$ rows, $n-k$ columns, and weight $c$. 
Let 
$M(k,n,c)$ be the set of all permutations $\pi \in S_n$
with $k$ weak excedences and $c$ crossings. 
Then there is a bijection $\Phi: T(k,n,c) \rightarrow M(k,n,c)$.
Moreover,  
$\Phi$  
maps a tableau $\T$ whose path $\phi(\T)$ has south steps in positions 
$S \subset \{1, \dots , n\}$ to a permutation 
with weak excedence set $S$.
\end{theorem}

Given $\tau = (\ttt_1, \dots , \ttt_n) \in \{0,1\}^{n}$,
we define $W(\tau)$ to be the subset of $\{1, 2, \dots , n+1\}$ which
contains $1$ and also contains $i+1$ if and only if $\tau_i = 1$.
Define the {\it excedence weight} $\wt(\pi)$ to be the number of crossings 
of $\pi$, and   
define $F'_{\tau}(q): = \sum_{\pi} q^{\wt(\pi)}$, 
where the sum is over all permutations in $S_{n+1}$ with 
weak excedence set $W(\tau)$.

By combining Corollary \ref{Cor1} and Theorem \ref{Bij}, we obtain
the following result.

\begin{corollary}\label{perminterp}
Fix $\ttt = (\ttt_1, \dots , \ttt_n) \in \{0,1\}^n$  and suppose
that $\alpha = \beta = 1$.
The probability of finding the PASEP model in configuration 
$(\ttt_1, \dots , \ttt_n)$ in the steady state is 
\begin{equation*}
\frac{F'_{\tau}(q)}{Z_n}.
\end{equation*}
Here $F'_{\tau}(q)$ is the generating function which enumerates 
permutations in $S_{n+1}$ with weak excedence set $W(\tau)$ according to number of
crossings.
\end{corollary}

This corollary can also be obtained in another way.
Namely, in \cite{Corteel2}, the authors provide a lattice path interpretation
for the probability of find the PASEP model in a particular configuration.
If one then applies a certain bijection on lattice paths 
\cite[Lemma 9]{Corteel} and then the 
bijection of Foata and Zeilberger \cite{FZ} from lattice paths to 
permutations, one arrives 
at Corollary \ref{perminterp} in this manner.

\begin{example}
Suppose that we are interested in $\ttt=(0,1,0)$.  If 
we set $\alpha = \beta = 1$ and use the results of Example
\ref{example1}, then we find that the probability of finding the PASEP
model in configuration $\ttt$ in the steady state (for $\alpha = \beta = 1$)
is 
$\frac{q+2}{2q^2+8q+14}$.
Note that $W(\tau) = (1,3)$.
Now observe that there are precisely three permutations in $S_4$
with weak excedence set $(1,3)$: $(3,1,4,2)$, and 
$(4,1,3,2)$, and $(2,1,4,3)$.  The first of these permutations
has one crossing and the second two have none, so the generating function
for permutations of type $\ttt$ is $q+2$, in agreement with
Corollary \ref{perminterp}.
\end{example}

We now introduce some more definitions concerning permutations.
If $\pi \in S_n$ is a permutation, we say that $\pi$ has a 
{\it descent} in position $i$ if $\pi(i) > \pi(i+1)$. 
The {\it descent set} of $\pi$ is the subset of $\{1, 2, \dots , n-1\}$
where $\pi$ has descents. 
And we say that $\pi$ has an occurrence of a {\it generalized pattern of type} 
$2-31$ if there is a pair $(i,j)$ such that $1 \leq i < j < n$ and  
$\pi(j+1)<\pi(i)<\pi(j)$.  These patterns were defined by 
Babson and Steingr\'{\i}msson \cite{BS}.

Given $\tau = (\ttt_1, \dots , \ttt_n) \in \{0,1\}^{n-1}$, 
we define $D(\tau)$ to be the subset of $\{1, 2, \dots , n \}$ which 
contains $i$ if and only if $\tau_i = 1$.
We define 
the {\it descent weight} $\wt'(\pi)$ to be the number of crossings 
of $\pi$, and   
define $F''_{\tau}(q): = \sum_{\pi} q^{\wt'(\pi)}$, 
where the sum is over all permutations in $S_{n+1}$ with
descent set $D(\tau)$.

We need the following result from \cite{SW}, which 
is also proved in
proved in \cite[Proposition 6]{Corteel}.

\begin{theorem}\cite[Theorem 16]{SW}.
There is a bijection $\Psi$ on permutations in $S_n$, which 
sends a permutation $\pi$ with descent set $D(\tau)$ to a permutation
$\Psi(\pi)$ with weak excedence set $W(\tau)$.  Moreover, the 
excedence weight of $\Psi(\pi)$ is equal to the descent weight of 
$\Psi(\pi)$, i.e.\ the number of crossings of $\Psi(\pi)$ is equal 
to the number of occurrences of the pattern $2-31$ in $\pi$.
\end{theorem}

By applying this result to Corollary \ref{perminterp}, we get the 
following.

\begin{corollary}\label{perminterp2}
Fix $\ttt = (\ttt_1, \dots , \ttt_n) \in \{0,1\}^n$ 
and suppose $\alpha = \beta = 1$.
The probability of finding the PASEP model in configuration 
$(\ttt_1, \dots , \ttt_n)$ in the steady state is 
\begin{equation*}
\frac{F''_{\tau}(q)}{Z_n}.
\end{equation*}
Here $F''_{\tau}(q)$ is the generating function which enumerates 
permutations in $S_{n+1}$ with descent set $D(\tau)$ 
according to the number of
occurrences of the pattern $2-31$.
\end{corollary}

Now we will show how these results imply a result of 
the first author \cite{Corteel}.

First we recall the definition of the polynomials
$\hat{E}_{k,n}(q)$ which were introduced by the second author \cite{Williams}. 
Define  
\begin{equation*}
{\hat{E}}_{k,n}(q) =
          q^{k-k^2} \sum_{i=0}^{k-1} (-1)^i [k-i]^n q^{ki-k}
       \left( {n \choose i} q^{k-i} + {n \choose i-1}\right).
\end{equation*} 

It was shown (implicitly) there (and more explicitly 
in \cite{SW}) that 
$\hat{E}_{k,n}(q)$ enumerates permutation tableaux  with 
$k$ rows and $n-k$ columns
according to weight.   Additionally, it was shown 
in \cite{Williams} that at $q=-1, 0, 1$, 
$\hat{E}_{k,n}(q)$ specializes to binomial coefficients, 
Naryana numbers, and Eulerian numbers. 

In \cite{Corteel}, the following connection was made between
these polynomials and the PASEP model.

\begin{theorem} \cite{Corteel} \label{Sylvie}
Let $\alpha = \beta = 1$.
Then in the steady state,
the probability that the PASEP model with 
$n$ sites is in a configuration with precisely $k$ particles is:
\begin{equation*}
\frac{\hat{E}_{k+1,n+1}(q)}{Z_n}
\end{equation*} 
\end{theorem}

Since the polynomials $\hat{E}_{k,n}(q)$ enumerate permutation 
tableaux with $k$ rows according to weight, 
Corollary \ref{Cor1} implies Theorem \ref{Sylvie}.

\section{Proof of Theorem \ref{Th1}}\label{Proof}

Recall the following requirement for permutation tableaux.
When we read the columns of a permutation tableau $\T$ from right to left,
if in any column  
we have a $0$ which lies beneath some $1$, then all 
entries to the left of $0$ (which are in the same row) must also
be $0$'s.

Also recall that an entry in a column of a permutation tableau is 
{\it restricted} if 
that entry is a $0$ which lies below some $1$. 
And a row is 
{\it unrestricted}
if it does not contain a restricted entry.   

\begin{proof}
We will prove Theorem \ref{Th1} inductively, by finding
a precise combinatorial interpretation 
for each of the entries in the top row of a matrix product such as
\begin{equation*}
\prod_{i=1}^n (\ttt_i D_1 + (1-\ttt_i) E_1).
\end{equation*}

Let 
$F^i_{\lambda}(q):= \sum_{\T} q^{\wt(\T)} \alpha^{-f(\T)} \beta^{-u(\T)}
= \sum_{\T} q^{\wt(\T)} \alpha^{-f(\T)} \beta^{-i+1}$, 
where the 
sum ranges over all permutation tableaux $\T$ of shape
$\lambda$ which have precisely $i$ unrestricted rows.  
Clearly $F_{\lambda}(q) =\sum_{i\geq 1}
F^i_{\lambda}(q)$.
Note that a permutation tableau will always have at least one 
unrestricted row (the top one).

Now let $M_{\lambda}$ be the matrix
$\prod_{i=1}^n (\ttt_i D_1 + (1-\ttt_i) E_1)$, where
$(\ttt_1, \dots , \ttt_n) = \phi(\lambda)$. 
We claim that the entry $M_{\lambda}[1,i]$ in position $(1,i)$ of $M_{\lambda}$
is $F^i_{\lambda}(q)$.

First note that this claim holds when $M_{\lambda}$ is equal to $D_1$
or $E_1$.  $M_{\lambda}$ is equal to $D_1$ when $\lambda$ is the 
partition $(0,0)$.  In this case there is a single  permutation
tableaux with two rows and no columns, which has weight $0$, 
and two unrestricted rows (hence $f(\T)=0$ and $u(\T)=1$).
This corresponds to the fact that the top row of $D_1$ is 
$(0,\beta^{-1},0,\dots)$.  
Similarly, $M_{\lambda}$ is equal to $E_1$ when $\lambda$ is the partition
$(1)$.  
In this case there is a single permutation tableaux with one
box which is filled with a $1$ (hence $f(\T)=1$).  
This tableaux has weight $0$ and 
one unrestricted row (hence $u(\T)=0$), 
corresponding to the fact that the top row
of $E_1$ is $(\alpha^{-1},0,0,\dots)$.

Using induction, assume that our claim is true for $\lambda$ with expanse 
less than or equal to $n$.  In other words, we can 
interpret the $i$th entry of the top row of $M_{\lambda}$ as 
a generating function enumerating permutation tableaux of shape
$\lambda$ with $i$
unrestricted rows, according to weight.
Let us now consider how these generating functions
will change if we instead consider permutation tableaux of 
shape ${\lambda}'$, where ${\lambda}'$ is a new shape obtained from 
$\lambda$ by either adding a new row of length $0$, or a new column 
whose length is the number of rows (including rows of length $0$) of  
$\lambda$.  The corresponding operation on paths is the following:
we take the path corresponding
to $\lambda$ and add an additional step from its south-west corner
which is either south or west.

It is easy to see how adding a step south to the partition path ---
i.e.\ adding an empty row to $\lambda$ --- will affect the generating
functions.  Any permutation tableau $\T$ of the new shape $\lambda'$ will
be a permutation tableau of shape
$\lambda$ with one additional unrestricted row
(the last row).  Therefore if $\lambda'$ is equal to $\lambda$ union
a row of length $0$, we will have 
$F_{\lambda'}^{1}(q)=0$, and $F_{\lambda'}^{i+1}(q)=
\beta^{-1} F_{\lambda}^{i}(q)$ for
$i\geq 1$.
This corresponds to the fact that in the matrix product
$M_{\lambda}D_1$, the new top row will be 
$(0,\beta^{-1} M_{\lambda}[1,1], \beta^{-1} M_{\lambda}[1,2], \dots )$.

It now remains to see how adding a step west to the partition path --- 
i.e.\ adding an extra column  of length $r$ to the left-hand-side of 
$\lambda=(\lambda_1, \dots , \lambda_r)$ 
--- will affect the generating functions.  
In this case our new partition $\lambda'$ is equal to 
$(\lambda_1 +1, \dots , 
\lambda_r +1)$.

Let $h_{a,b}(q)$ denote the polynomial 
$\beta^{b-a} \left(\alpha^{-1} q^{a-1} {b-1 \choose a-1} +
\sum_{j=0}^{a-2} q^j {b-a+j \choose j} \right)$ for $a \leq b$.
We claim that 
$F_{\lambda'}^{a}(q) = \sum_{b \geq a} h_{a,b}(q) F_{\lambda}^{b}(q)$.

To prove this, it is enough to show the following.
Fix a permutation tableau $\T$ of shape $\lambda$ which has precisely
$b$ unrestricted rows.  Consider all ways of adding an additional
(maximal) column $C$ to the left of $\T$, in order to build a new permutation
tableau $\T'$ of shape $\lambda'$ which has precisely $a$ unrestricted
rows.  Then the generating function for these new tableaux according 
to weight is precisely $h_{a,b}(q)q^{\wt(\T)}$.

To prove this last statement, consider the process of adding the column
$C$ to $\T$. Since $\T$ has $b$ unrestricted rows, there are only $b$ entries
in $C$ in which we can choose to put either a $1$ or $0$; all other
entries are forced to be $0$.  Let us number those unrestricted entries from top
to bottom by 
$c_1$ to $c_b$.  Since we want our new column $C$ to 
add an additional $b-a$ restricted positions, in our filling 
of the entries $c_1,\dots,c_b$, we must have precisely $b-a$ $0$'s below
a $1$.  If the top-most $1$ in our filling of $c_1, \dots ,c_b$ is 
in position $c_i$, then the $i-1$ entries above it are all $0$'s, and
we have precisely ${b-i \choose b-a}$ ways to choose which entries to make
$0$ below  $c_i$.  The other $(b-i)-(b-a)=a-i$ entries must be $1$'s.
This particular choice of column will therefore contribute the extra weight 
$q^{a-i}$ to the weight of $\T$ -- if $i \neq 1$ -- and 
will contribute the weight $\alpha^{-1} q^{a-1}$ -- if $i=1$.
If we sum over all possible columns
$C$ which we may add to $\T$, we get the following:
$\alpha^{-1} q^{a-1} {b-1 \choose b-a} + 
\sum_{i=2}^a q^{a-i} {b-i \choose b-a}$ which is equal to 
$\alpha^{-1} q^{a-1} {b-1 \choose a-1} + 
\sum_{j=0}^{a-2} q^j {b-a+j \choose j}$.
This completes the proof that 
$F_{\lambda'}^{a}(q) = \sum_{b \geq a} h_{a,b}(q) F_{\lambda}^{b}(q)$.
And now note that this corresponds to the fact that in the matrix product
$M_{\lambda}E_1$, the  entry in the first row and $a$th column 
will be 
$\sum_{b \geq a} h_{a,b}(q) M_{\lambda}[1,b]$.

This now completes our proof that for any partition $\lambda$,
the entry $M_{\lambda}[1,i]$ in position $(1,i)$ of $M_{\lambda}$
is $F^i_{\lambda}(q)$.  And since $F_{\lambda}(q) = \sum_{i \geq 1}
F^i_{\lambda}(q)$, the theorem follows.
\end{proof}

\section{Applications of permutation tableaux}\label{Applications}

In this section we will give some applications of 
the connection of permutation tableaux 
to the PASEP.  The first application is a partial order 
on states of the PASEP and some monotonicity results (with respect to that 
partial order) of probabilities
of observing these states.  Note that some of the results below hold for
general $\alpha, \beta$ ($0 < \alpha \leq 1$ and $0 <\beta \leq 1$), and for others we need to assume that 
$\alpha=\beta=1$.

\begin{definition}
Let $\tau, \tau' \in \{0,1\}^n$ be two states of the PASEP which 
contain exactly $k$ particles.  We define the 
partial order $\prec$ by $\tau \prec \tau'$
if and only if $\lambda(\tau) \subset \lambda(\tau')$.  
\end{definition}

\begin{figure}[h]
  \centering
\includegraphics[height=4in]{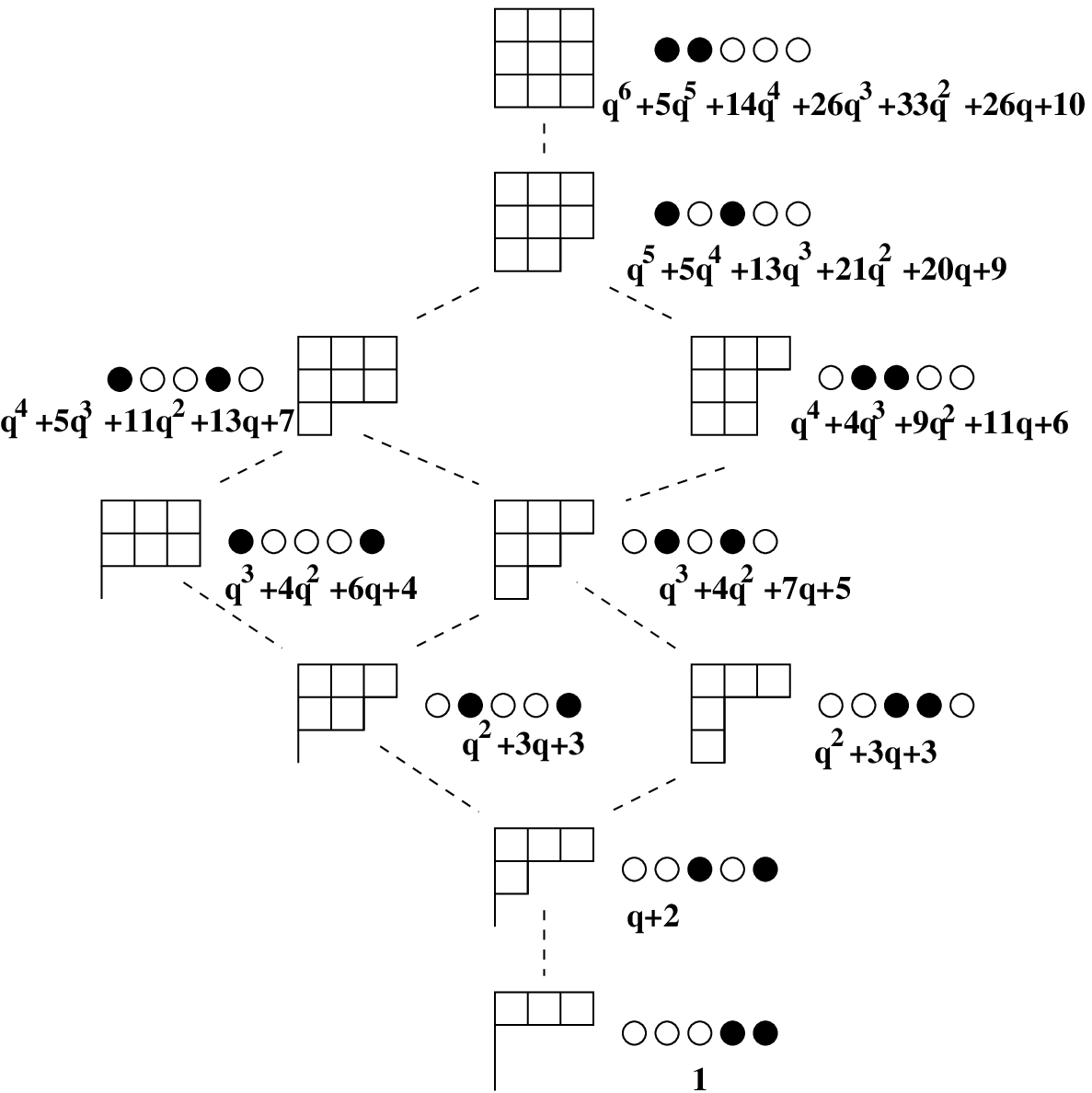}
\caption{A partial order for states of the PASEP}
\label{Order}
\end{figure}

Figure \ref{Order} illustrates this partial order when $n=5$ and
$k=2$, and it shows the (unnormalized) probabilities
$f_n(\tau)$ that each of these states occurs, for $\alpha=\beta=1$.

We now give two simple inequalities relating 
$f_n(\tau)$ to $f_n(\tau')$ when $\tau \prec \tau'$.

\begin{proposition}\label{qdiff}
Let $\alpha$ and $\beta$ be general.  Suppose that $\tau \prec \tau'$ and let 
$d:= |\lambda(\tau')| - |\lambda(\tau)|$.  That is, 
$d$ is the difference between the cardinalities of 
the Young diagrams $\lambda(\tau')$ and $\lambda(\tau)$.
Then $f_n(\tau') - q^d f_n(\tau)$ is a non-negative polynomial 
in $q$, $\alpha^{-1}$, and $\beta^{-1}$.
\end{proposition}

\begin{proof}
Observe that any permutation tableau of shape 
$\tau$ can be naturally extended to a permutation tableau
of shape $\tau'$ by filling in all boxes of 
$\lambda(\tau') \setminus \lambda(\tau)$ with $1$'s.  
This new permutation tableau has weight $q^d$ times the weight
of the old one.
\end{proof}

\begin{proposition}\label{mono}
Let $\alpha=\beta=1$.
Suppose that $\tau \prec \tau'$.  Then 
$f_n(\tau') - f_n(\tau)$ is a non-negative polynomial in $q$.
In other words, as one moves up the partial order
$\prec$, the coefficients of $f_n(\tau)$ monotonically
increase.
\end{proposition} 

Proposition \ref{mono} is false for general $\alpha$ and $\beta$.
For example, $f_4(1,1,0,0)-f_4(1,0,1,0)$ is a Laurent polynomial
in $q$, $\alpha$, and $\beta$ with some negative coefficients.

We will defer the proof of Proposition \ref{mono} to the next
section because it is most easily proved using Motzkin paths.

We thank the referee for pointing out that the corollary below
follows from our previous results.

\begin{corollary}
Let $\alpha=\beta=1$.  Suppose that $\tau \prec \tau'$.  Let 
$d$ be any integer such that 
$0 \leq d \leq |\lambda(\tau')| - |\lambda(\tau)|$.  
Then
$f_n(\tau')-q^d f_n(\tau)$ is a non-negative polynomial in $q$.
\end{corollary}

\begin{proof}
This follows from Propositions \ref{qdiff} and \ref{mono}.  Simply
choose a saturated chain of partitions 
$\tau = \tau_0 \prec \tau_1 \prec \dots \tau_m=\tau'$
such that $|\lambda(\tau_{k+1})|-|\lambda(\tau_k)|=1$
and apply either Proposition \ref{qdiff} or Proposition 
\ref{mono} at each step. 
\end{proof}

The above results  make sense intuitively,
since in a model with particles entering from the left and leaving
to the right, it is more likely that a given particle will be 
further to the left than to the right.  Also recall that 
the probability of hopping left ($\frac{q}{n+1}$ for some 
$q\leq 1$) is at most the 
probability of hopping right ($\frac{1}{n+1}$)

Another application of permutation tableaux
is that these objects
allow one to 
read off the main
recurrences for the PASEP 
(\cite[Theorem 1]{Corteel2}) quite easily.
The rest of the results in this section hold for general 
$\alpha$ and $\beta$.

\begin{figure}[h]
  \centerline{\epsfig{figure=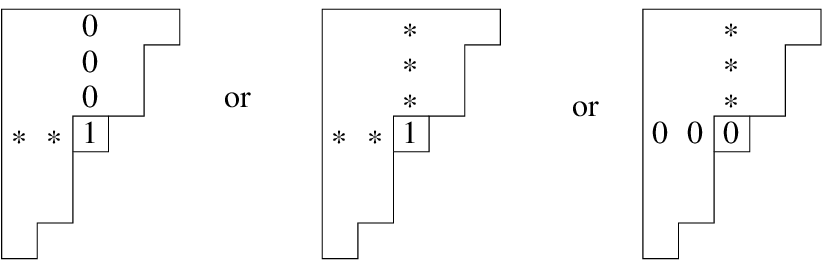}}
\caption{Recurrence for $F_{\lambda}(q)$}
\label{Recurrence}
\end{figure}

Fix a partition $\lambda = (\lambda_1, \dots , \lambda_k)$, 
and let $\tau = (\tau_1, \dots , \tau_n)$ be the 
unique vector in $\{0,1\}^n$
such that $\lambda = \lambda (\tau)$.
Choose any {\it corner} of $\lambda$, i.e.\ 
the last box of a row $\lambda_i$ in $\lambda$
such that $\lambda_{i+1} < \lambda_i$.
(Equivalently, a pair $(j,j+1)$ of entries in $\tau$ such that
$\tau_j = 1$ and $\tau_{j+1}=0$.)
As Figure \ref{Recurrence} illustrates, there is a simple
recurrence for $F_{\lambda}(q)$.

Explicitly, any valid filling of $\lambda$ is obtained in one
of the following ways:
\begin{itemize}
\item inserting a column whose bottom entry is $1$ and whose
other entries are $0$ after the $(\lambda_i - 1)$st column of
a valid filling of $(\lambda_1 - 1, \lambda_2 - 1, \dots , 
\lambda_k-1)$;
\item adding a $1$ to the end of the $i$th row of a valid filling
of the shape $(\lambda_1, \lambda_2, \dots , \lambda_i - 1, \dots , \lambda_k)$;
\item inserting an all-zero row of length $\lambda_i$ as the $i$th row
in a valid filling of $(\lambda_1, \dots , \lambda_{i-1}, \lambda_{i+1}, 
\dots , \lambda_k)$.
\end{itemize}

Note that none of the operations above change the number of $1$'s in 
the first row or the number of unrestricted rows.  Therefore, 
for $1 < i \leq k$ and $\lambda_i > \lambda_{i+1}$, we have that 
\begin{equation*}
F_{\lambda}(q) = F_{(\lambda_1 - 1, \dots , \lambda_i -1, \lambda_{i+1}, \dots , \lambda_k)} + 
F_{(\lambda_1, \dots , \lambda_{i-1}, \lambda_{i+1}, \dots , \lambda_k)} +
qF_{(\lambda_1, \dots , \lambda_{i-1},\lambda_i -1, \lambda_{i+1}, \dots,
\lambda_k)}.
\end{equation*}

By Corollary \ref{Cor1}, this recurrence translates into a
recurrence for the PASEP model.
Recall our notation $\frac{f_n(\tau_1, \dots , \tau_n)}{Z_n}$ for 
the probability of finding the system in 
configuration $(\tau_1, \dots , \tau_n)$ in the steady state.

\begin{corollary}\cite[Theorem1]{Corteel2}
Let $\alpha$ and $\beta$ be general.  Then
\begin{align*}
f_n(\tau_1, \tau_2, \dots , \tau_{j-1}, 1, 0, 
\tau_{j+2}, \dots , \tau_n) = 
&f_{n-1}(\tau_1, \tau_2, \dots , \tau_{j-1}, 1,
\tau_{j+2}, \dots , \tau_n) +\\
&qf_n(\tau_1, \tau_2, \dots , \tau_{j-1},0,1,
\tau_{j+2},\dots , \tau_n) +\\
&f_{n-1}(\tau_1, \dots , \tau_{j-1}, 0, 
\tau_{j+2}, \dots , \tau_n).
\end{align*}
\end{corollary}

It is also easy to prove the following, by consideration 
of permutation tableaux.

\begin{corollary}\cite[Theorem 1]{Corteel2}
Let $\alpha$ and $\beta$ be general.  We have that 
$f_n(\tau_1,  \dots , \tau_{n-1},1)=
\frac{1}{\beta} f_n(\tau_1,  \dots , \tau_{n-1})$
and also 
$f_n(0, \tau_1, \tau_2, \dots , \tau_n)=
\frac{1}{\alpha} f_n(\tau_1, \tau_2, \dots , \tau_n)$.
\end{corollary}

\begin{proof}
This also follows from
Corollary \ref{Cor1}.  The first equality holds because
the  generating function for permutation tableaux of shape
$(\lambda_1, \lambda_2, \dots , \lambda_{n-1},0)$ is clearly equal 
to the  generating function for permutation tableaux of shape
$(\lambda_1, \lambda_2, \dots , \lambda_{n-1})$ except that
there is one additional unrestricted row.  The second
equality holds because  
the generating function for permutation tableaux of shape
$(\lambda_1 + 1, \lambda_2, \dots , \lambda_{n})$ is equal 
to the generating function for permutation tableaux of shape
$(\lambda_1, \lambda_2, \dots , \lambda_{n})$ except that there
is one additional entry in the first row which is a $1$. 
(Recall that each column of a permutation tableaux is required to 
contain exactly one $1$ and that the weight of a column is equal
to the number of $1$'s it contains {\it beyond} that requisite
one.)
\end{proof}

\section{Connection with Motzkin paths}\label{Motzkin}

We now give another solution to the matrix ansatz, 
which is essentially the one 
described in
Derrida et al \cite{Derrida1}, and which 
has an interpretation in terms of 
Motzkin paths.  

Note that in this section we restrict to the 
case $\alpha = \beta = 1$.  These results can be extended to the more
general case, but we restrict ourselves to this case since our main
purpose here is to obtain the result needed for the
monotonicity result of the previous section.

Let $D_0$ be the (infinite) upper triangular
matrix $(d_{ij})$ such that 
$d_{i,i}=[i]$ where $[i]$ is the $q$-analog of $i$, 
$d_{i,i+1}=[i+1]$, and $d_{i,j}=0$ for $j\neq i, i+1$.

That is, $D_0$ is the matrix
\[ \left( \begin{array}{ccccc}
[1] & [2] & 0 & 0 & \dots \\
0 & [2] & [3] & 0 & \dots \\
0 & 0 & [3] & [4] & \dots \\
0 & 0 & 0 & [4] & \dots \\
\vdots & \vdots & \vdots & \vdots &
\end{array} \right). \]

Let $E_0$ be the (infinite) lower triangular matrix
$(e_{i,j})$ such that 
$e_{i,i}=[i]$, 
$e_{i+1,i}=[i]$, and $e_{i,j}=0$ for $j\neq i,i-1$.

That is, $E_0$ is the matrix
\[ \left( \begin{array}{ccccc}
1 & 0 & 0 & 0 & \dots \\
1 & [2] & 0 & 0 & \dots \\
0 & [2] & [3] & 0 & \dots \\
0 & 0 & [3] & [4] & \dots \\
\vdots & \vdots & \vdots & \vdots &
\end{array} \right). \]

Let $W_0$ be the (row) vector $(1,0,0,\dots )$ and $V_0$ be the (column)
vector
$(1,0,0,\dots)$.
It is now easy to check that 
$D_0 E_0 -q E_0 D_0 = D_0 + E_0$, 
$DV_0 = V_0$, and $W_0 E = W_0$.

We will now give an interpretation of the steady states
of the PASEP
in terms of bicolored
Motzkin paths.  This result is very similar to a result obtained in 
Brak et al \cite{Corteel2}, via a combinatorial derivation.  
Note also that in \cite{BrakEssam}, Brak and Essam
considered the case $q=0$ of the PASEP and 
gave multiple interpretations for the steady states  
in terms of various weighted lattice paths.

We define a  {\it bicolored Motzkin path} of length $n$ to be a sequence of 
steps in the plane
$c=(c_1, \dots , c_n)$ such that 
$c_i \in \{N, S, E, \bar{E}\}$ for $1 \leq i \leq n$, 
which starts and ends at height $0$, and 
always stays at or above height $0$.  
That is, if we define the 
height at the $i$th step to be 
$h_i = \{j < i \ \vert \ c_j = N \} - \{j < i \ \vert \ c_j=S\}$,
then $h_i \geq 0$ for $1 \leq i \leq n$, and furthermore, 
$h_{n+1}=0$.

We now assign a weight $\wt(p)$ to a bicolored Motzkin path $p$ as
follows: step $c_i$ is assigned weight $[h+1]$ if $c_i$ ends at height
$h$, and the weight of $c$ is defined to be the product of the weights of
all the $c_i$'s.

Fix a sequence $\ttt=(\ttt_1, \dots , \ttt_n) \in \{0,1\}^n$.
We say that a bicolored Motzkin path $c$ has {\it type} 
$\ttt$ if $c_i=N$ or $E$ whenever $\ttt_i=1$ and 
$c_i=S$ or $\bar{E}$ whenever $\ttt_i=0$. 

We now define $F'''_{\tau}(q)$ to be the generating function for all 
bicolored Motzkin paths of type $\tau$: that is, 
$F'''_{\tau}(q) :=\sum_p q^{\wt(p)}$ where the sum is over all bicolored
Motzkin paths of type $\tau$.

The main result of this section is the following.

\begin{proposition}\label{M}
Choose a sequence $\ttt=(\ttt_1, \dots , \ttt_n) \in \{0,1\}^n$.
Then 
\begin{equation*}
F'''_{\tau}(q) = W_0 (\prod_{i=1}^n (\ttt_i D_0 + (1-\ttt_i) E_0)) V_0
\end{equation*}
\end{proposition}

Let us make the convention in this section that rows and columns 
of our matrices are indexed by non-negative integers (including $0$).
Then Proposition \ref{M} says that  
the generating function $F'''_{\tau}(q)$ is equal 
to the entry in the $0$th row and $0$th column of the matrix
product 
$\prod_{i=1}^n (\ttt_i D_0 + (1-\ttt_i) E_0)$.  

\begin{proof}
Using the definition of matrix multiplication, we
expand $\prod_{i=1}^n (\ttt_i D_0 + (1-\ttt_i) E_0)$
as the sum of terms of the form
$A_1[1,i_1] A_2[i_1, i_2] \dots A_n[i_{n-1},1]$
where the $A_j$'s are matrices $D_0$ or $E_0$, depending on
the sequence $\tau$.  It is obvious that these terms correspond
to the Motzkin paths of type $\tau$.
\end{proof}

As before, we can  use 
Proposition \ref{M}
together with Theorem \ref{ansatz} to deduce the following.

\begin{corollary}\label{Cor3}
Fix $\ttt = (\ttt_1, \dots , \ttt_n) \in \{0,1\}^n$. 
The probability of finding the PASEP model in configuration 
$(\ttt_1, \dots , \ttt_n)$ in the steady state is 
\begin{equation*}
\frac{F'''_{\tau}(q)}{Z_n}.
\end{equation*}
Here $F'''_{\tau}(q)$ is the generating function for 
bicolored Motzkin paths of type $\tau$.
\end{corollary}

We now use this result to prove Proposition \ref{mono}
from the previous section:

\begin{proof}
It is sufficient to consider two configurations
$\tau$ and $\tau'$ whose corresponding diagrams
$\lambda(\tau)$ and $\lambda(\tau')$ differ by one box.
That is, there exists some $i$ such that:\\
$\bullet$ $\tau_i = 0$, $\tau_{i+1}=1$\\
$\bullet$ $\tau'_i = 1$, $\tau'_{i+1}=0$\\
$\bullet$ $\tau'_j = \tau_j$ for $j \neq i, i+1$.

We now observe that any bicolored Motzkin path $p$ of 
type $\tau$ (with weight $w$) can be mapped 
to a bicolored Motzkin path $p'$
of type $\tau'$ whose weight $w'$ is coefficient-wise
greater than $w$.

If $p=(c_1, \dots , c_n)$ then we let 
$p':=(c_1, \dots, c_{i-1},c_{i+1},c_{i},c_{i+2},\dots,c_n)$.
That is, $p'$ is the path obtained from $p$ by switching the 
$i$th and $i+1$st steps.
It is clear that whether we switch $SN$ steps for $NS$,
or $SE$ for $ES$, or $\bar{E}N$ for $N\bar{E}$, or 
$\bar{E}E$ for $E\bar{E}$, the resulting path is a valid
Motzkin path.  The weight $w'$ will differ from $w$ by
(respectively): replacing a factor 
$[h-1][h]$ with $[h+1][h]$; 
replacing a factor $[h-1]^2$ with $[h][h-1]$;
replacing a factor $[h][h+1]$ with $[h+1]^2$;
nothing.
Using the fact that $[n+1]-[n]$ is nonnegative, all of the
differences
$[h+1][h]-[h-1][h]$, 
$[h][h-1]-[h-1]^2$, $[h+1]^2-[h][h+1]$
are nonnegative; therefore in all cases,
$w'-w$ will be a nonnegative polynomial.
By Corollary \ref{Cor3}, we are done.
\end{proof}

Note that one can use the same argument to give another proof of
Proposition \ref{qdiff} in the case that $\alpha=\beta=1$. We repeat the
same argument but now 
$w'$ differs from $w$ by: 
replacing $q[h-1][h]$ with $[h+1][h]$; replacing 
$q[h-1]^2$
with $[h][h-1]$; and replacing $q[h][h+1]$ with $[h+1]^2$. Then all
differences of such polynomials will be nonnegative, since $[n+1]-q[n]$ is
nonnegative.

\begin{remark}
We now summarize our combinatorial interpretations for the 
steady state probabilities of the PASEP model; these are given by
our solutions $(D_j, E_j)$ to the matrix 
ansatz (where $\alpha = \beta = 1$), for $j=0,1$.
Note that since 
\begin{equation*}
W_j (\prod_{i=1}^n (\ttt_i D_j + (1-\ttt_i) E_j)) V_j 
\end{equation*}
describes the steady state probability of the PASEP model
for both $j=0,1$, these two formulas must be equal.  
In particular, both of these are formulas for all  of the following:
\begin{itemize}
\item $F_{\lambda(\tau)}(q)$, the weight-generating function for 
permutation tableaux of shape $\lambda(\tau)$.
\item $F'_{\tau}(q)$, the generating function for permutations of 
excedence type $\tau$, enumerated according to crossings.
\item $F''_{\tau}(q)$, the generating function for permutations of 
descent type $\tau$, enumerated according to occurrences of $2-31$.
\item $F'''_{\tau}(q)$, the generating function for weighted
bicolored Motzkin 
paths.
\end{itemize}

Moreover, the product $W_j (D_j+E_j)^n V_j$ is the weight-generating 
function for all of the following:
\begin{itemize}
\item permutation tableaux with expanse $n+1$.
\item permutations in $S_{n+1}$, enumerated according to crossings. 
\item permutations in $S_{n+1}$, enumerated according to occurrences of 
$2-31$.
\item weighted bicolored Motzkin paths of length $n$.
\end{itemize}

\end{remark}

\end{document}